\documentclass[amsmath,amssymb,superscriptaddress,longbibliography,reprint,10pt]{revtex4-1}
\usepackage[T1]{fontenc}
\usepackage[utf8]{inputenc}
\usepackage{helvet}

\pdfoutput=1
\usepackage[dvipsnames,rgb,dvips]{xcolor}
\usepackage{graphicx}
\usepackage{bbold}
\usepackage{float}
\usepackage{overpic}
\usepackage{amsmath,amssymb,amsfonts,MnSymbol}
\usepackage{hyperref}
\makeatletter
\def\amsbb{\use@mathgroup \M@U \symAMSb}
\makeatother
\usepackage{mathrsfs}
\usepackage{overpic}

\begin{document}
\title{\large Extra-factorial sum: a graph‐theoretic \\ parameter in Hamiltonian cycles of complete weighted graphs}
\author{V. Papadinas}
\affiliation{Department of Informatics, Hellenic Open University, Patras, Greece}
\author{W. Xiong}
\affiliation{Institute for Theoretical Physics, Heidelberg University, Heidelberg, Germany}
\author{N. A. Valous}
\affiliation{National Center for Tumor Diseases, German Cancer Research Center, Heidelberg, Germany}

\begin{abstract}
A graph‐theoretic parameter, in a form of a function, called the extra-factorial sum is discussed. The main results are presented in ref. \cite{RefJMain} (Nastou \textit{et al.}, Optim Lett, 10, 1203--1220, 2016) and the reader is strongly advised to study the aforementioned paper. The current work presents subject matter in a tutorial form with proofs and some newer unpublished results towards the end (lemma six extension and lemma seven). The extra-factorial sum is relevant to Hamiltonian cycles of complete weighted graphs $WH_n$ with $n$ vertices and is obtained for each edge of $WH_n$. If this sum is multiplied by $1 / (n - 2)$ then it gives directly the arithmetic mean of the sum of lengths $l_i$ of all Hamiltonian cycles that traverse a selected edge $e_q$. The number of terms in this sum is a factorial proven to be $(n - 2)!$ which signifies that its value depends on $n$. Using the extra-factorial sum, the arithmetic mean of the sum of the squared lengths of $(n - 1)! / 2$ Hamiltonian cycles of $WH_{n}$ can be obtained as well. 
\\
\\
\scriptsize{Authors to whom all correspondence should be addressed; e-mails: farmermath70@gmail.com\textsuperscript{1} and nek.valous@nct-heidelberg.de\textsuperscript{3}}
\end{abstract}

\maketitle
\begin{figure}[h]
    \centering
    \includegraphics[width = 0.47\textwidth]{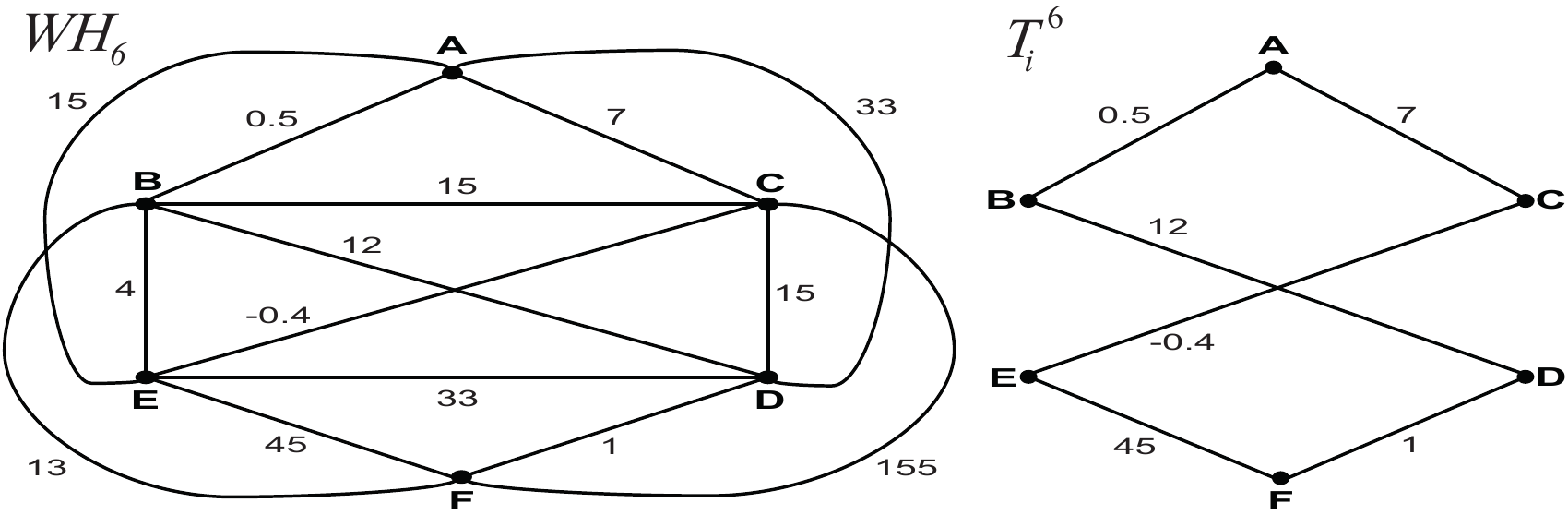}
	\caption{Complete weighted graph $WH_6$ and Hamiltonian cycle $T_i^6$ (vertices $=$ edges) originating from $WH_6$.}
	\label{fig1}
\end{figure}

Consider a complete graph $H_n$ with $n$ vertices. Such a graph has always $n\times (n - 1) / 2$ edges, since for each pair of vertices $(x, y)$ with $x \neq y$, exists only one edge $e_q$ with $1 \leq q \leq n \times (n - 1) / 2$ that connects them. Additionally, each $e_q$ can have a weight $w(e_q) \in$ $\mathbb{R}$ which makes $H_n$ a weighted graph $WH_n$. Fig. \ref{fig1} shows a $WH_6$. In a graph $H_n$, every closed walk $T_i$ with $1 \leq i \leq (n - 1)! / 2$ that traverses each vertex exactly once starting and ending at the same vertex is called a Hamiltonian cycle \cite{RefB1}. Each $T_i$ is comprised of $e_T$ edges and $n$ vertices $(e_T = n)$. The $i$th $T_i$ with $n$ vertices is called $T_i^n$. The length $l_i$ of a $T_i^n$ is obtained by summing the relevant $w(e_q)$. Fig. \ref{fig1} also shows a $T_i^6$: \textbf{ACEFDBA} which originates from $WH_6$ and has $l_i = 65.1$. The number of $T_i^n$ is given by $(n - 1)! / 2$ \cite{RefB1}, e.g. $WH_6$ has 60 $T_i^6$.

A graph‐theoretic problem is introduced which is defined as follows: for selected $e_q$ in $WH_n$, the aim is to find the sum of lengths $l_i$ of the $T_i^n$ that traverse $e_q$. Furthermore, it will be proved that each $e_q$ has a number of $T_i^n$ that traverse it which is $(n - 2)!$. The Hamiltonian path problem is NP-complete \cite{RefJ1} and summing the $l_i$ of the $(n - 2)!$ $T_i^n$ that traverse $e_q$ is computationally expensive. On the other hand, the arithmetic mean of this sum can be obtained and this new function is defined as the extra-factorial sum \cite{RefDef}. This is a fraction with numerator the sum of $l_i$ of the $(n - 2)!$ $T_i^n$ traversing $e_q$ and denominator the value $(n - 3)!$. The parameter can be obtained directly without the need of computing the numerator and denominator separately. It will be proven that the extra-factorial sum of any $e_q$, if multiplied by $1 / (n - 2)$, yields the arithmetic mean of $l_i$ of the $(n - 2)!$ $T_i^n$ traversing $e_q$. A direct deduction is that for any $WH_n$, the extra-factorial sum can be obtained for each $e_q$. This graph‐theoretic parameter can be presented in a 2D Cartesian chart with the $y$-axis showing the values of the extra-factorial sum in the interval $(-\infty, +\infty)$ and the $x$-axis showing the ranked $e_q$ based on the extra-factorial sum values. Such a curve is a visualization of the distribution of $l_i$ in relation to the edges, and is predominantly a qualitative measure for comparing different $WH_n$.

\textbf{Lemma 1}. \textit{A Hamiltonian cycle generator GT is a $T_i^n$ in which if a vertex is added then this cycle can create new $T_i^{n+1}$ with $(n + 1)$ edges 
and $(n + 1)$ vertices. After creating the new $T_i^{n+1}$, GT is replaced by its child cycles.}

\textbf{Proof}. 
If vertex X is added to $T_i^n$ then for each $e_q$ there is a unique pair ($e_{q1}$, $e_{q2}$) with common X that can break the edges, thus creating the new $T_i^{n+1}$. 
In this way, $n$ new $T_i^{n+1}$ can be created with $(n + 1)$ edges and $(n + 1)$ vertices each. $\blacksquare$

\begin{figure}[h]
    \centering
    \includegraphics[width = 0.45\textwidth]{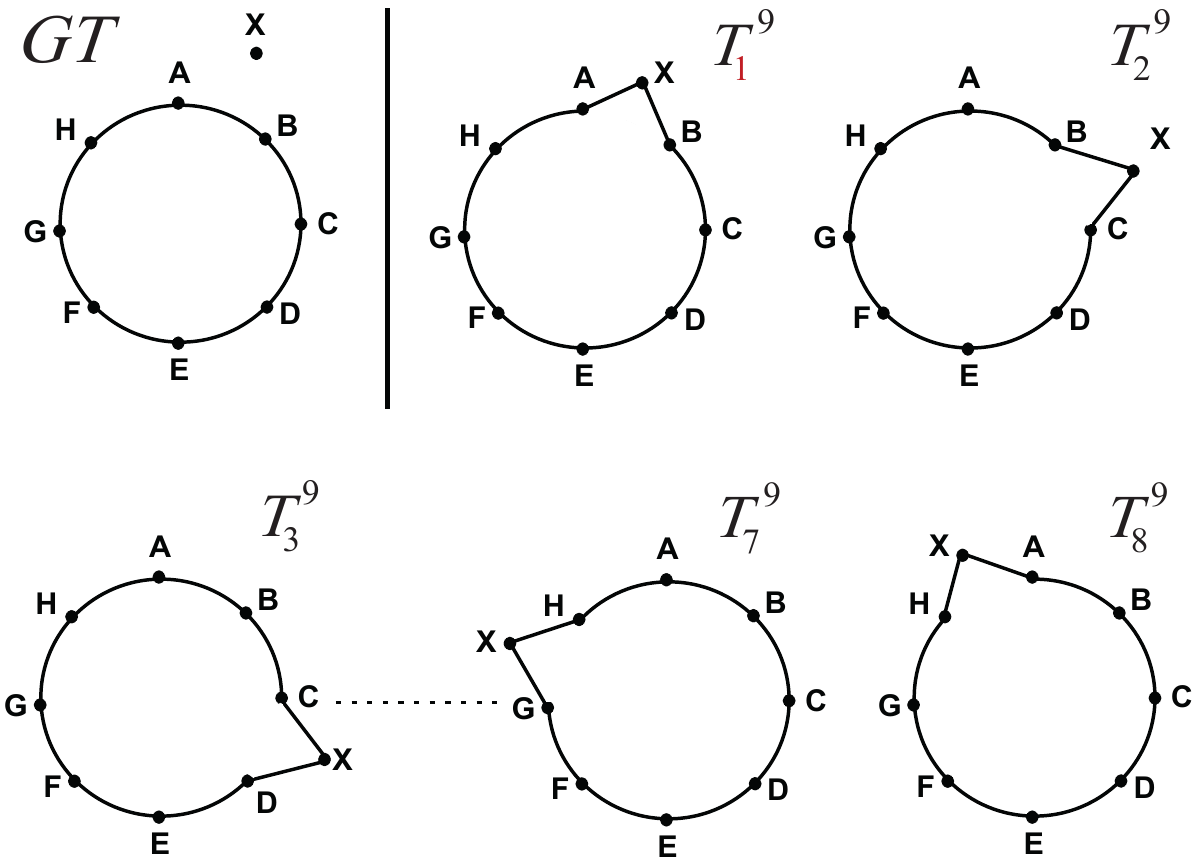}
	\caption{Creation of eight new $T_i^{8 + 1}$ starting with an initial cycle generator $T_i^8$ and vertex X outside the cycle.}
	\label{fig2}
\end{figure}

With an initial cycle generator $T_i^8 =$ \textbf{ABCDEFGHA} and vertex X outside the cycle, then new $T_i^{8 + 1}$: $T_1^{8 + 1}$, $T_2^{8 + 1}$, $T_3^{8 + 1}$,$T_4^{8 + 1}$, $T_5^{8 + 1}$, $T_6^{8 + 1}$, $T_7^{8 + 1}$, $T_8^{8 + 1}$ are created (Fig. \ref{fig2}). Each new $T_i^{n+1}$ traverses the vertices of GT including X. In essence, each $e_q$ of GT (\textbf{AB, BC, CD, DE, EF, FG, GH, HA}) breaks in order to create new $T_i^{n+1}$. For example, the creation of \textbf{AXBCDEFGHA} requires the breaking of \textbf{AB} and the addition of \textbf{AX} and \textbf{XB}. The initial \textbf{ABCDEFGHA} ceases to exist after the creation of the new $T_i^{n+1}$. The obtained $T_i^{8 + 1}$ are equal to the $e_q$ of $T_i^8$, hence $T_i^8$ is defined as a GT. The process of creating $T_i^{n+1}$ is described in Table \ref{tab1} with the Steps Insert Vertex Algorithm (SIVA). The inputs are $T_i^{n}$ and vertex X outside the cycle and the output is the $n$ new $T_i^{n+1}$ with $(n + 1)$ edges and $(n + 1)$ vertices \cite{RefJ2}.

\textbf{Corollary 1}. 
\textit{Every $T_i^n$ belongs to generation $(n - 3)$ defined as $G = (n - 3)$. The $T_i^n$ belonging to the same generation have the same number of edges and vertices. Since every $T_i^n$ has same number of edges and vertices then minimum is  $n = 3$. $T_i^3$ cannot be a GT therefore the generation is defined as $G = (n - 3) = 0$. The $T_i^n$ with $G = 0$ is the initial $T_1^3$.}

\begin{table}[h]
\caption{The process of creating $T_i^{n+1}$ with the Steps Insert Vertex Algorithm (SIVA).}
\label{tab1}
\begin{tabular*}{\columnwidth}{@{\extracolsep{\fill}}p{\dimexpr 0.03\linewidth-2\tabcolsep}p{\dimexpr 0.97\linewidth-2\tabcolsep}}
\\ \hline
\multicolumn{2}{p{\dimexpr 0.97\linewidth-2\tabcolsep}}{\textbf{Input}: $T_i^n$ $(n \geq 3)$ and vertex X outside the cycle.}\\
$\rightarrow$ & \textbf{Step 1}: copy $T_i^n$ and vertex X ($n$ times).\\
$\rightarrow$ & \textbf{Step 2}: for each copy $T_i^n$ (step 1), remove a different $e_q$.\\
$\rightarrow$ & \textbf{Step 3}: for each $T_i^n$ (step 2), add a copy of X.\\
$\rightarrow$ & \textbf{Step 4}: for each $i$th copy (step 3), add a new pair $(e_{q1}, e_{q2})$ so that each broken edge is relinked with X.\\
\multicolumn{2}{p{\dimexpr 0.975\linewidth-2\tabcolsep}}{\textbf{Output}: $n$ new $T_i^{n+1}$ $(1 \leq i \leq n)$ with $(n + 1)$ edges and $(n + 1)$ vertices.}
\\ \hline
\end{tabular*}
\end{table}

Henceforth, the edges of $H_n$ or $WH_n$ are denoted by $e_n$ and the edges of $T_i^n$ by $e_T^n$. When SIVA is applied to $T_i^n$ ($n \geq 3$) then this is defined as the $GT$ with $G = (n - 3)$ that can create $n$ new $T_i^{n+1}$ with $G + 1$ having $(n + 1)$ vertices and $(n + 1)$ edges. SIVA can be extended by creating $(n - 1)! / 2$ $T_i^n$ that exist in $H_n$ with $n$ vertices and $n (n - 1) / 2$ edges $(n \geq 3)$. Using SIVA to count $T_i^n$ offers the opportunity to count only the $T_i^n$ for which certain categorization criteria may apply. For $H_3$, the relationship $[n (n - 1) / 2] = e_n = n = 3$ is valid; this means that $H_3$ has the additional property of being a $T_1^3$ with $G = 0$, since for $T_i^n$ the number of vertices is equal to the number of edges.

\textbf{Lemma 2}. \textit{For $H_3$ and $k$ vertices $(k \geq 1)$ outside the graph, $H_3$ is a $T_1^3$ (initial $GT$). The $[(3 + k) - 1]! / 2$ $T_i^n$ of the new $H_{k + 3}$ are counted one-by-one when the new $k$ vertices are added to $H_3$ iteratively using SIVA.}

\begin{figure}[h]
    \centering
    \includegraphics[width = 0.40\textwidth]{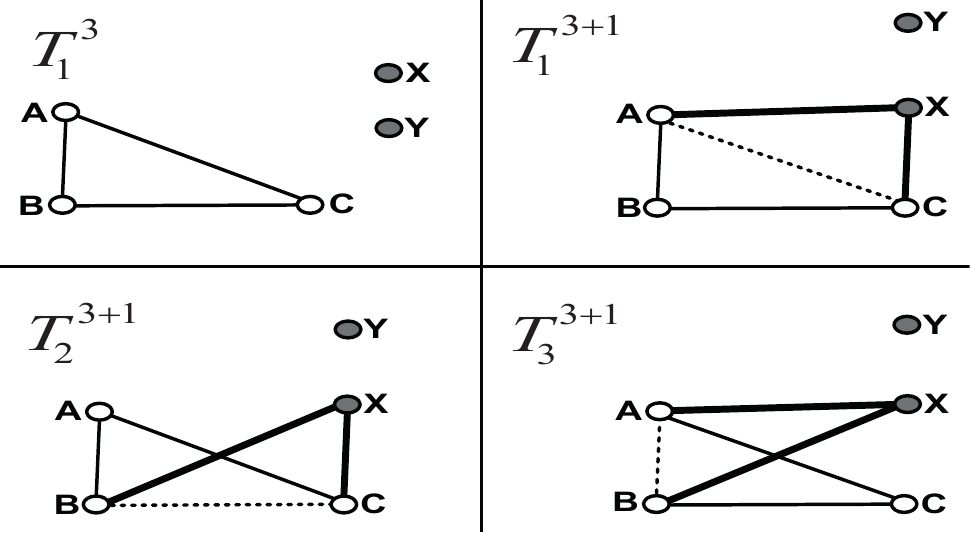}
	\caption{Creation of three new $T_i^4$ ($G = 1$).}
	\label{fig3}
\end{figure}

\textbf{Proof}. A graph $H_3$ with $n = e_n$ is also a $T_1^3$. SIVA starts from $T_1^3$ and iterates $k$ times to introduce, in each repetition, a new vertex in every $T_i^{3 + j}$ with $G = [(3 + j) - 3]$, for $1 \leq j \leq k$. Fig. \ref{fig3} shows $H_3 = \textbf{ABC}$ and two new vertices X and Y $(k = 2)$; \textbf{ABC} is also a $T_1^3 = \textbf{ABCA}$. 

\begin{figure}[h]
    \centering
    \includegraphics[width = 0.45\textwidth]{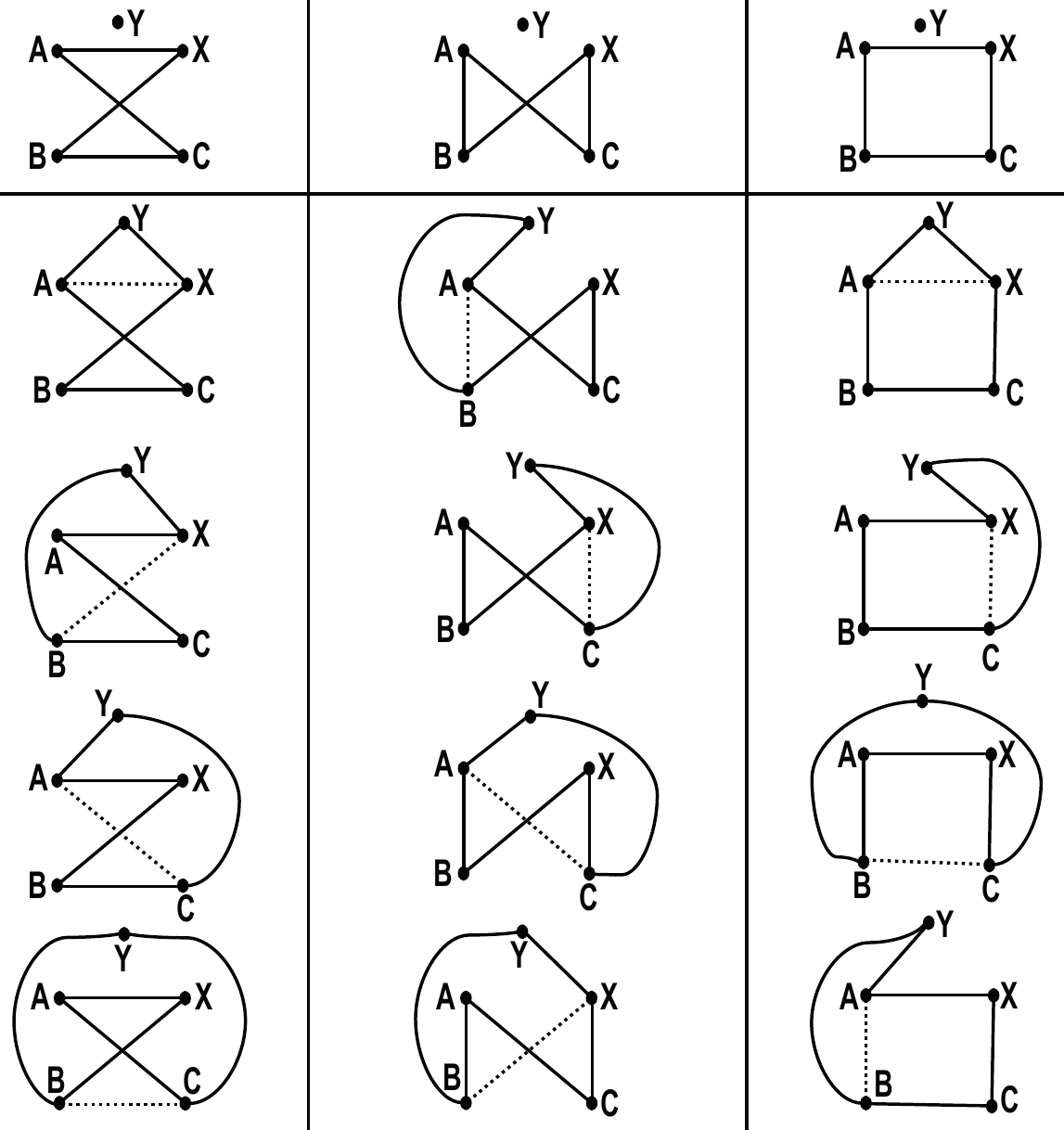}
	\caption{Creation of twelve new $T_i^5$ ($G = 2$).}
	\label{fig4}
\end{figure}

\noindent
Applying SIVA results in adding X into \textbf{ABCA}, which leads to the creation of three new $T_i^4$ ($G = 1$): $T_1^4 = \textbf{AXCBA}$, $T_2^4 = \textbf{ACXBA}$, and $T_3^4 = \textbf{AXBCA}$. Dotted edges represent breaking locations for creating the new cycles. After inserting X, \textbf{ABC} is converted to $H_4 = \textbf{ABCX}$. Continuing along, SIVA adds Y to each $T_i^4$ $(1 \leq i \leq 3)$ of the first generation \textbf{ABCX}. After adding Y, twelve new $T_i^5$ ($G = 2$) are created ($1 \leq i \leq 12$) (Fig. \ref{fig4}). 

\begin{figure}[h]
    \centering
    \includegraphics[width = 0.35\textwidth]{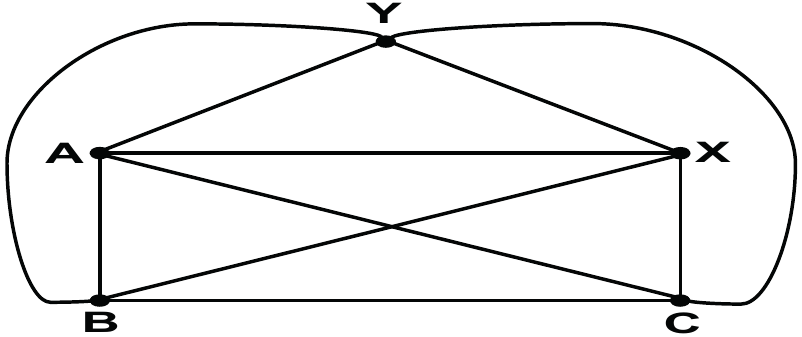}
	\caption{The newly created $H_5 =$ \textbf{ABCXY}.}
	\label{fig5}
\end{figure}

\noindent
The new $T_i^5$ are created when SIVA is applied to each $T_i^4$ of \textbf{ABCX}. After inserting Y, $H_4$ is converted to $H_5 = \textbf{ABCXY}$ (Fig. \ref{fig5}). In the general case (starting with $n = 3$), when adding a new vertex to $H_n$ then SIVA places the vertex to every $T_i^n$, for $1 \leq i \leq [(n - 1)! / 2]$. After adding the $k$th vertex, new $H_{n+k}$ are created with $(n + k)$ vertices and $[(n + k) (n + k - 1)] / 2$ edges, for $(n + k) > 3$ and $k \geq 1$. This new graph incorporates the set of $[(n + k) - 1]! / 2$ $T_i^{n + k}$. Each $T_i^{n + k}$ belongs to $G = (n + k - 3)$ with $(n + k)$ vertices and $(n + k)$ edges. Hence, the following recursive function is obtained:

\begin{equation}
 \frac{[(n-1)-1]!}{2}e_{T^{n-1}} = \frac{(n-1)!}{2}
\end{equation}

Alternatively, this can be expressed as:

\begin{equation}
 \frac{(n_{old}-1)!}{2}e_{T^{n_{old}}} = \frac{(n_{new}-1)!}{2}
\end{equation}

\noindent
where $n_{new} - n_{old} = 1$. This algebraic relationship states that the product of $(n_{old} - 1)! / 2$ of the $T_i^{n_{old}}$ by $e_T^{n_{old}}$ in each $T_i^{n_{old}}$ yields $(n_{new} - 1)! / 2$ of the $T_i^{n_{new}}$; these are created if $H_n$ increases its vertices by $1$ and its edges by $n$. Validity can be proved algebraically as a function of $n$, when SIVA increases the vertices of $H_n$ by $k$. The new $T_i^5$ ($G = 2$) as a function of $n$, for $n = 3$, are $(1 \times n) \times (n + 1) = 12$. The term $(1 \times n)$ are the $T_i^4$ ($G = 1$) and the term $(n + 1)$ are the $e_T^4$ in each $T_i^4$ ($G = 1$). With the addition of vertices and starting with $n = 3$, the new $T_i^{n+k}$ tends exactly to:

\begin{equation}
n(n+1)(n+2)(n+3)\cdots (n+k-1) = \frac{(n+k-1)!}{2}
\end{equation}

\noindent
This concludes the proof. $\blacksquare$

\textbf{Lemma 3}. \textit{The number of $T_i^n$ in $H_n$ $(n \geq 3)$ which traverse selected $e_q$ is $(n - 2)!$.}

\begin{figure}[h]
    \centering
    \includegraphics[width = 0.42\textwidth]{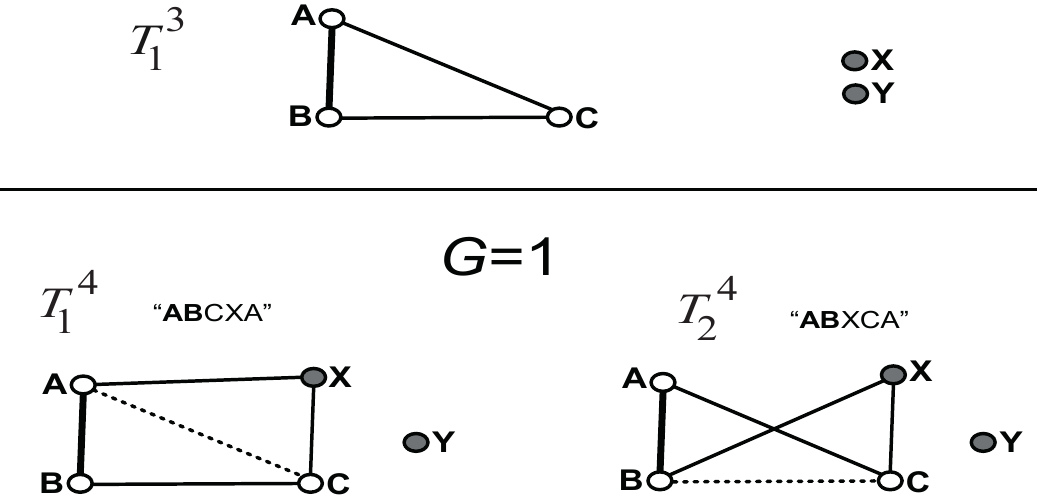}
	\caption{Creation of new $T_i^4$ ($G = 2$) that traverse selected $e_q^{AB}$.}
	\label{fig6}
\end{figure}

\begin{figure}[h]
    \centering
    \includegraphics[width = 0.42\textwidth]{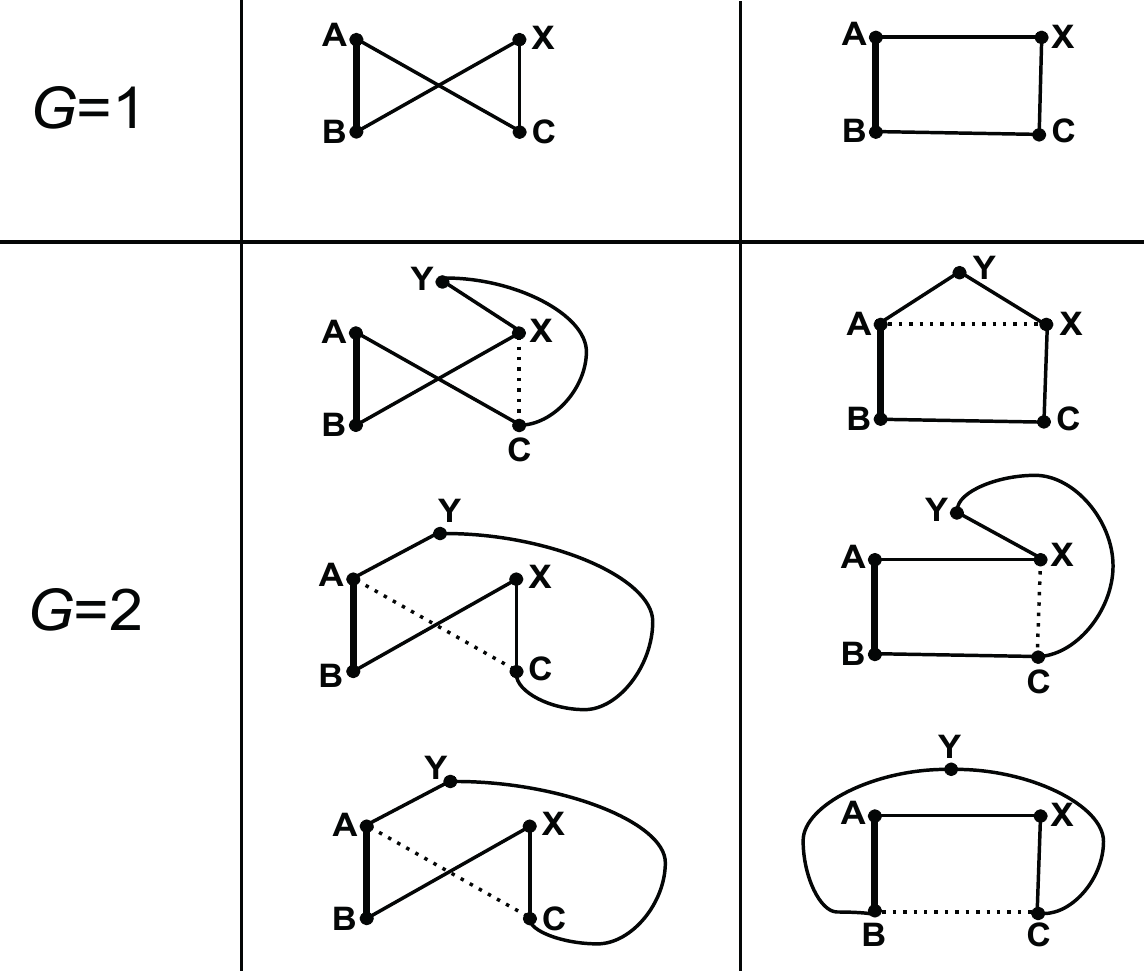}
	\caption{Creation of new $T_i^5$ ($G = 2$) that traverse selected $e_q^{AB}$.}
	\label{fig7}
\end{figure}

\textbf{Proof}. Counting the $T_i^n (n \geq 3)$ can be done using SIVA but with some differences comparing to Lemma 2: $e_q$ does not break. Specifically, if $e_q$ is selected then $k = n - 3$ and the first $T_1^3$ is created with the two vertices of $e_q$ and any one vertex of $H_n$. Then SIVA iterates by inserting the remaining $k$ vertices for $n = (3 + k)$. In each iteration, it is imperative that $e_q$ does not break. This variation is visualized in Fig. \ref{fig6} and Fig. \ref{fig7} for $G = 1$ and $G = 2$, respectively. For $H_5$ (Fig. \ref{fig5}) and selected $e_q^{AB}$ that links \textbf{A} and \textbf{B}, in order to count the $T_i^n$ that always traverse $e_q^{AB}$, only the $T_i^{n-1}$ traversing $e_q^{AB}$ are counted recursively. The procedure begins with $H_5$ and $e_q^{AB}$; then three of the five vertices are selected such that $e_q^{AB}$ is always traversed: these three vertices form an $H_3$ and a $T_1^3$. SIVA is applied to the edges ($e_q^{AC}$, $e_q^{BC}$) of $T_1^3$ except $e_q^{AB}$; hence a $T_1^4$ ($G = 1$) is created. For $T_1^4 = \textbf{ABCXA}$ and $T_2^4 = \textbf{ABXCA}$, SIVA adds the fifth and final vertex of $H_5$ without breaking $e_q^{AB}$. In the general case of $H_n$, the algorithm proceeds till all remaining vertices are added. Hence, a set of $T_i^n$ is created that always traverse $e_q^{AB}$. The number of $T_i^n$ is $(n - 2)!$ which is proven by induction: this is valid for $H_3$ since the cycles that traverse $e_q$ are $(n - 2)! = 1$. Similarly for $H_4$ since $(n - 2)! = 2$. Therefore:

\begin{equation}
[(n-1)-2]!(e_{T^{n-1}}-1) = (n-2)!
\end{equation}

\noindent
The term $(n - 2)!$ expressing the $T_i^n$ that traverse $e_q$ can be obtained using SIVA, for each $T_i^{n-1}$ traversing $e_q$ without breaking it. The fact that $e_q$ does not break is denoted by the second term $(e_{T^{n-1}}-1)$: this expresses the edges in $T_i^{n-1}$ allowed to break in an one-by-one fashion for creating the new $T_i^n$. Overall, the product defines the $T_i^{n-1}$ traversing $e_q$ by the value $(e_{T^{n-1}}-1)$ from each $T_i^{n-1}$ excluding $e_q$. Since $e_{T^n} = n$, this can be written as:

\begin{equation}
[(n-1)-2]!(e_{T^{n-1}}-1) = (n-3)!(n-2) = (n-2)!
\end{equation}

Alternatively, this can be expressed as:

\begin{equation}
 (n_{old}-2)!(e_{T^{n_{old}}}-1) = (n_{new}-2)!
 \end{equation}

\noindent
where $n_{new} - n_{old} = 1$. This expresses the $T_i^{n-1}$ in $H_{n-1}$ that always traverse $e_q$ by the edges in $T_i^{n-1}$ excluding $e_q$; this equals to the $T_i^n$ in $H_n$ that traverse $e_q$. $\blacksquare$

An edge intersects or not any other edge; if $e_q^{AB}$ intersects \textbf{A} and \textbf{B} (Fig. \ref{fig5}) then it intersects an edge adjacent to $e_q^{AB}$, meaning that they share the same vertex. Hence, for every edge intersecting $e_q^{AB}$: {$e_q^{AX}$, $e_q^{BY}$} where \textbf{$X, Y \neq  (A \parallel B)$}, e.g. edges intersecting \textbf{AB} are: \textbf{AY, AX, AC} and \textbf{BY, BX, BC}.

\textbf{Lemma 4}. \textit{The number of $T_i^n$ in $H_n$ $(n \geq 3)$ which traverse a selected pair of adjacent edges is $(n - 3)!$.}

\textbf{Proof}.
Counting the $T_i^n$ ($n \geq 3$) that satisfy the lemma is proven through Lemma 3. For example, in $H_5$ (Fig. \ref{fig5}) edges $e_q^{AY}$ and $e_q^{YX}$ with common \textbf{Y} are selected. Initially, it is observed that there are no $T_i^5$ that traverse $e_q^{AX}$ and no edge capable of linking \textbf{Y} with vertices other than \textbf{A} and \textbf{X}, therefore $e_q^{AX}$, $e_q^{YB}$, and $e_q^{YC}$ have to be removed. For $H_n$ $(n \geq 4)$ and pair of adjacent edges (e.g. $e_q^{AY}$ and $e_q^{YX}$), the edge that links the vertices of the selected adjacent edges is removed. This means $e_q^{AX}$ and any other edge that links \textbf{Y} with any other vertex excluding $e_q^{YK}$, where \textbf{$K \neq (A \parallel B)$}. After the removal of edges, \textbf{Y} can be removed temporarily since for all  $T_i^5$, \textbf{A} and \textbf{B} are always linked. The graph resulting from the removal of \textbf{Y} is an $H_4$ with $(n - 1) = 4$ vertices. In this new $H_4$, the number of $T_i^4$ that traverse $e_q^{AB}$ is equal to the number of $T_i^5$ traversing the adjacent $e_q^{AY}$ and $e_q^{YX}$ of the initial $H_5$. Previously, it was proven that the number of $T_i^n$ that traverse a selected edge is $(n - 2)!$. Consequently, the number of $T_i^n$ that traverse any selected pair of adjacent edges is $(n - 1) - 2! = (n - 3)!$. $\blacksquare$

\textbf{Lemma 5}. \textit{The number of $T_i^n$ in $H_n$ ($n \geq 4$) which traverse a selected pair of non-adjacent edges is $2 (n - 3)!$.}

\begin{figure}[h]
    \centering
    \includegraphics[width = 0.42\textwidth]{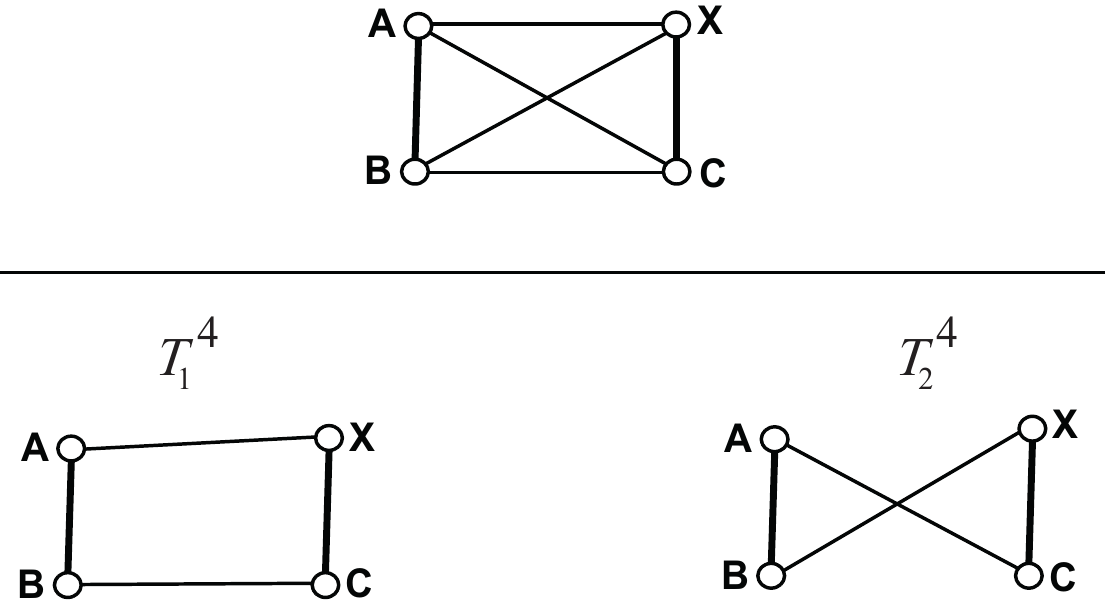}
	\caption{Creation of new $T_i^4$ that traverse non-adjacent $e_q^{AB}$ and $e_q^{XC}$.}
	\label{fig8}
\end{figure}

\begin{figure}[h]
    \centering
    \includegraphics[width = 0.42\textwidth]{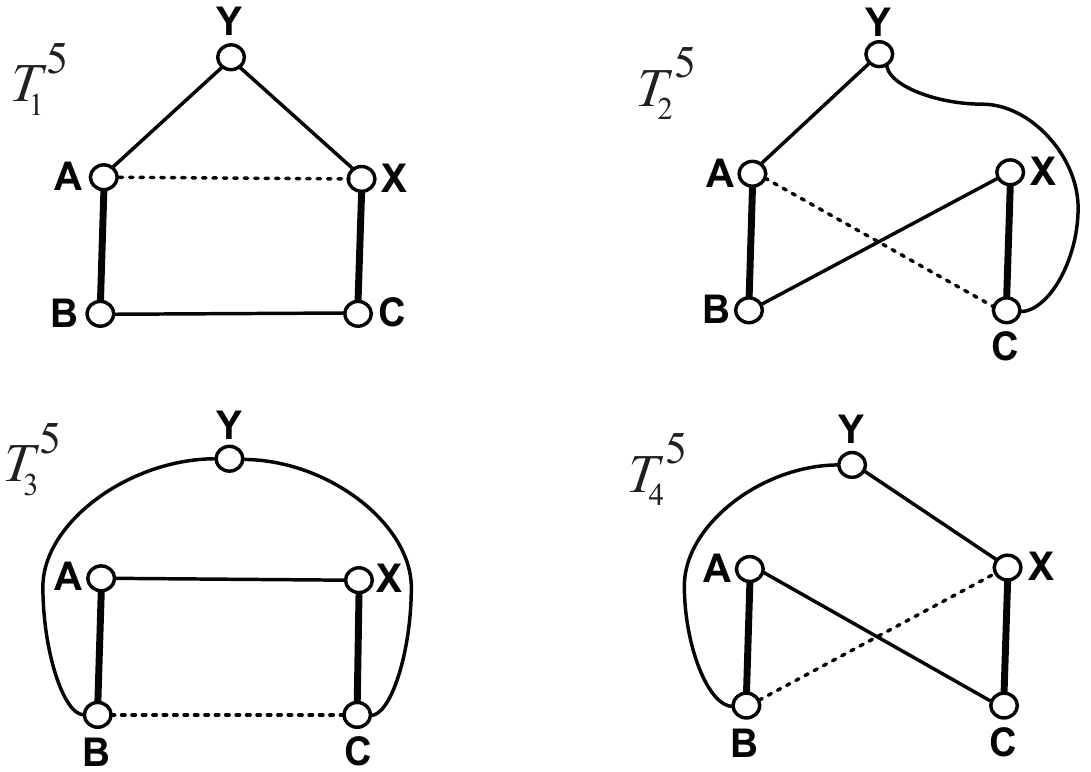}
	\caption{Creation of new $T_i^5$ that traverse non-adjacent $e_q^{AB}$ and $e_q^{XC}$.}
	\label{fig9}
\end{figure}

\textbf{Proof}. For $T_i^n$ ($n \geq 4$) and two selected non-adjacent edges, then SIVA creates $(e_T^n - 2) = (n - 2)$ new $T_i^{n+1}$ with $(n + 1)$ vertices since the two non-adjacent edges do not break. Counting the $T_i^n$ that fulfill the lemma starts from $n = 4$ since this is the minimum value of $n$ permitting the creation of $T_i^n$ that can have pairs of non-adjacent edges. Only $T_1^4$, $T_2^4$, and $T_3^4$ exist in $H_4$; from these and for every pair of non-adjacent edges (e.g. $e_q^{AB}$ and $e_q^{XC}$) only $T_1^4$ and $T_2^4$ traverse them, therefore counting starts from them (Fig. \ref{fig8}). SIVA creates the $T_i^n$ contained in each of the initial $T_1^4$ and $T_2^4$ without breaking $e_q^{AB}$ and $e_q^{XC}$. SIVA iterates $k$ times with $n = k - 4$ for each $T_1^4$ and $T_2^4$ and creates $2 [(4 + k) - 3]!$ new $T_i^{k+4}$ with $(4 + k)$ vertices and $(4 + k)$ edges. For $k = 1$ and $n = 4$, the expression $2 [(n + k) - 3]! = 4$ is valid since for $T_1^4 = \textbf{ABCDA}$ and $T_2^4 = \textbf{ABDCA}$ with the addition of a new vertex, SIVA creates $(e_T^n - 2) = (n - 2) = 2$ new $T_i^5$ that traverse $e_q^{AB}$ and $e_q^{XC}$. Given that $T_1^4$ and $T_2^4$ are two then $2 (e_T^n - 2) = 2 (n - 2) = 4$ new $T_i^5$ are created (Fig. \ref{fig9}). 

\noindent
The dotted lines (Fig. \ref{fig9}) denote the edges that break in order to add the new $Y$. The number of $T_i^5$ that fulfill the lemma is $2 (n - 3)! = 4$. When SIVA adds a new vertex in $H_5$ (without breaking $e_q^{AB}$ and $e_q^{XC}$) then for $n = 5$ and $k = 1$, $2 [(n + k) - 3]! = 12$ new $T_i^6$ are created traversing $e_q^{AB}$ and $e_q^{XC}$ since for each of the old $T_i^5$ three vertices break. This is expressed as: $4 (e_T^n - 2) = 4 (n - 2) = 12$, and in the general case:

\begin{equation}
\begin{split}
& 2[(n-1)-3]!(e_{T^{n-1}}-2) = \\
& 2[(n-1)-3]![(n-1)-2] = \\
& 2(n-4)!(n-3)= \\
& 2(n-3)!
\end{split}
\end{equation}

\noindent
The constant expresses the $T_i^4$ where counting starts from and proceeds with SIVA $k$ times without breaking the non-adjacent edges. The term $(n - 4)!$ expresses the child cycles originating from each initial $T_i^4$ in $H_{n - 1}$. The term $(n - 3)$ expresses the edges belonging to $T_i^{n -1}$ that break in an one-by-one fashion (except for the selected edges) in order for SIVA to insert a new vertex. $\blacksquare$

\begin{figure}[h]
    \centering
    \includegraphics[width = 0.45\textwidth]{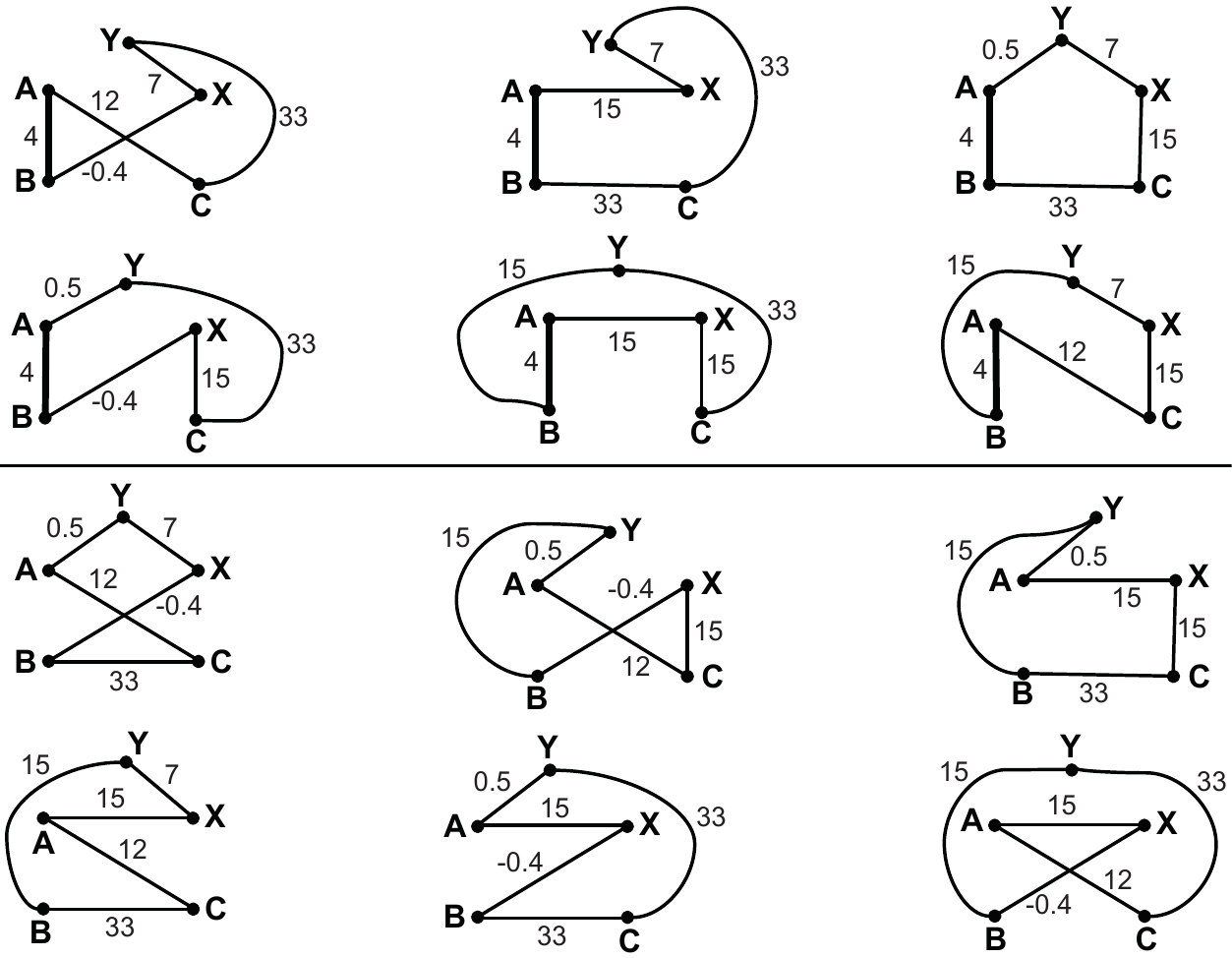}
	\caption{Twelve $T_i^5$ ($WH_5$) grouped into two categories.}
	\label{fig10}
\end{figure}

Consider the graph $WH_5$ (based on the graph $H_5$ of Fig. \ref{fig5}) with $w(e_q^{AB}) = 4$, $w(e_q^{AX}) = 15$, $w(e_q^{AC}) = 12$, $w(e_q^{AY}) = 0.5$, $w(e_q^{BC}) = 33$, $w(e_q^{BX}) = -0.4$, $w(e_q^{BY}) = 15$, $w(e_q^{CY}) = 33$, and $w(e_q^{XY}) = 7$. This graph has $n (n - 1) / 2 = 10$ edges and $(n - 1)! / 2 = 12$ $T_i^5$. Each $T_i^5$ has $l_i$ equal to the sum of its consisting edge weights. Fig. \ref{fig10} shows the $T_i^5$ into two groups; the first (top two rows) contains the cycles that traverse the selected $e_q^{AB}$ while the second (bottom two rows) the remaining that do not. 

\noindent
The number of $T_i^5$ that traverse $e_q^{AB}$ is $(n - 2)!$. According to Lemma 4, for the $T_i^5$ that traverse $e_q^{AB}$ the number of edges that intersect this edge is $(n - 3)!$ meaning that $e_q^{AY}$, $e_q^{AX}$, $e_q^{AC}$, $e_q^{BY}$, $e_q^{BX}$, and $e_q^{BC}$ appear twice in the $T_i^5$ subset that traverse $e_q^{AB}$. According to Lemma 5, for the $T_i^5$ that traverse $e_q^{AB}$ the number of edges that do not intersect this edge is $2 (n - 3)!$ meaning that $e_q^{YX}$, $e_q^{XC}$, $e_q^{YC}$ appear four times in the $T_i^5$ subset that traverse $e_q^{AB}$. The sum of $l_i$ of the $T_i^5$ that traverse $e_q^{AB}$ is: $\sum_{i=1}^{(n-2)!} l_{i}^{AB} = l_1 + l_2 +l_3 +l_4 +l_5 +l_6 = 349.2$, where $T_1^5 = \textbf{YCABXY} (l_1 = 55.6)$, $T_2^5 = \textbf{YCBAXY} (l_2 = 92)$, $T_3^5 = \textbf{YXCBAY} (l_3 = 59.5)$, $T_4^5 = \textbf{YCXBAY} (l_4 = 52.1)$, $T_5^5 = \textbf{YCXABY} (l_5 = 82)$, and $T_6^5 = \textbf{YXCABY} (l_6 = 53)$.

\textbf{Lemma 6}. \textit{For $WH_n$ and selected edge (e.g. $e_q^{AB}$) there exist a unique summational graph $\sum WH_n^{AB}$ that corresponds to $e_q^{AB}$. Each edge in $WH_n$ corresponds to a unique $\sum WH_n$. For $e_q^{AB}$, the resulting $\sum WH_n^{AB}$ is a copy of the initial $WH_n$ where each weight is multiplied as follows: i) $w(e_q^{AB})$ multiplied by $(n - 2)!$, ii) the weight of each edge intersecting $e_q^{AB}$ multiplied by $(n - 3)!$, and iii) the weight of each edge not intersecting $e_q^{AB}$ multiplied by $2 (n - 3)!$. The sum of $w(e_i)$ of $\sum WH_n$ corresponding to $e_q^{AB}$ is equal to the sum of $l_i$ of the $T_i^n$ that traverse $e_q^{AB}$:}

\begin{equation}
\sum_{i=1}^{n(n-1)/2} w(e_i) = \sum_{i=1}^{(n-2)!} l_i^{AB} 
\end{equation}

\begin{figure}[h]
    \centering
    \includegraphics[width = 0.42\textwidth]{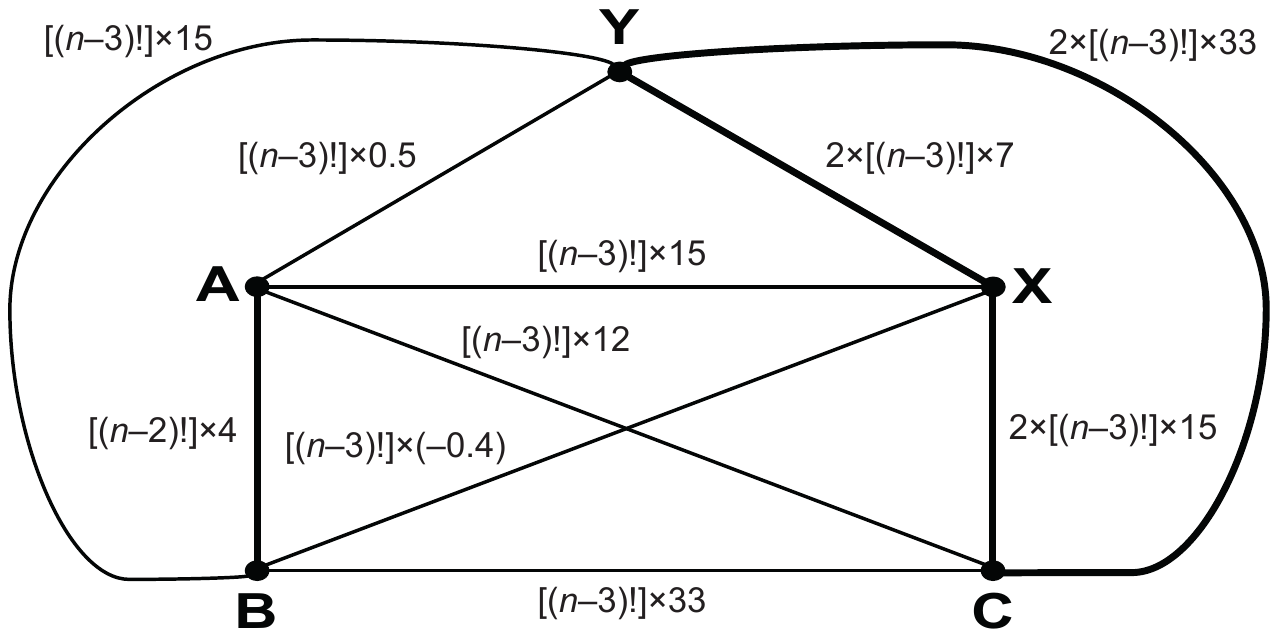}
	\caption{The summational graph $\sum WH_5^{AB}$ (\textbf{ABCXY}).}
	\label{fig11}
\end{figure}

\textbf{Proof}. The sum of $l_i$ of the $T_i^5$ traversing $e_q^{AB}$ is visualized in $\sum WH_5^{AB}$ (Fig. \ref{fig11}): this shows a copy of $WH_5$ for $e_q^{AB}$ where $w(e_q^{AB})$ is multiplied by $(n - 2)! = 6$, the weight of each edge intersecting $e_q^{AB}$ by $(n - 3)! = 2$, and the weight of the remaining edges by $2 (n - 3)! = 4$. The sum of edge weights in $\sum WH_5^{AB}$ is equal to the sum of $l_i$ of the $T_i^5$ that traverse $e_q^{AB}$.

\noindent
This concludes the proof. $\blacksquare$

The relationship between the sum of $l_i$ of the $T_i^n$ ($WH_n$) and the sum of edge weights in $\sum WH_n^{AB}$ is given by \cite{RefJ2}:

\begin{equation}
\begin{split}
\sum_{i=1}^{(n-2)!} l_{i}^{AB} = & \\
& (n-2)! X_1 + (n-3)! X_2 + 2(n-3)! X_3
\end{split}
\end{equation}

\noindent
where $X_1$ is the weight of $e_q^{AB}$, $X_2$ is the sum of weights with $X, Y \neq  (B \parallel A)$ for the edges that intersect $e_q^{AB}$ (their number is $2 (n - 2)$), and $X_3$ is the sum of weights with $X, Y \neq  (B \parallel A)$ for the edges that do not intersect $e_q^{AB}$ (their number is $[n (n - 1)] / 2 - [2 (n - 2)] - 1 = [(n - 2) (n - 3)] / 2$). Overall:

\begin{equation}
\begin{split}
X_1 = & w(e_q^{AB}) \\
X_2 = & \sum_{q_1=1}^{(n-2)} w(e_{q_1}^{AX})+ \sum_{q_2=1}^{(n-2)} w(e_{q_2}^{BY}) \\ 
X_3 = & \sum_{q=1}^{(n-2)(n-3)/2} w(e_q^{XY})
\end{split}
\end{equation}

\noindent
Since, the common term is $(n - 3)!$ the extra-factorial sum of edge \textbf{AB} can be written as:

\begin{equation}
\begin{split}
(!)\sum_{i=1}^{(n-2)!} l_{i}^{AB} & = (n-2)w(e_q^{AB}) + \\ 
& + [\sum_{q_1=1}^{(n-2)} w(e_{q_1}^{AX})+\sum_{q_2=1}^{(n-2)} w(e_{q_2}^{BY})] + \\
& + 2[\sum_{q=1}^{(n-2)(n-3)/2} w(e_q^{XY})]
\end{split}
\end{equation}

\noindent
The symbol (!) denotes the removal of the term $(n - 3)!$. The extra-factorial sum for a selected edge multiplied by $1 / (n - 2)$ is equal to the arithmetic mean of $l_i$ of the $(n - 2)!$ cycles that traverse that edge:

\begin{equation}
\begin{split}
(!)\sum_{i=1}^{(n-2)!} l_{i}^{AB} \frac{1}{(n-2)} & = \frac{\sum_{i=1}^{(n-2)!} l_{i}^{AB}}{(n-3)!} \frac{1}{(n-2)} = \\ 
& = \frac{\sum_{i=1}^{(n-2)!} l_{i}^{AB}}{(n-2)!}
\end{split}
\end{equation}

\begin{figure}[h]
    \centering
    \includegraphics[width = 0.33\textwidth]{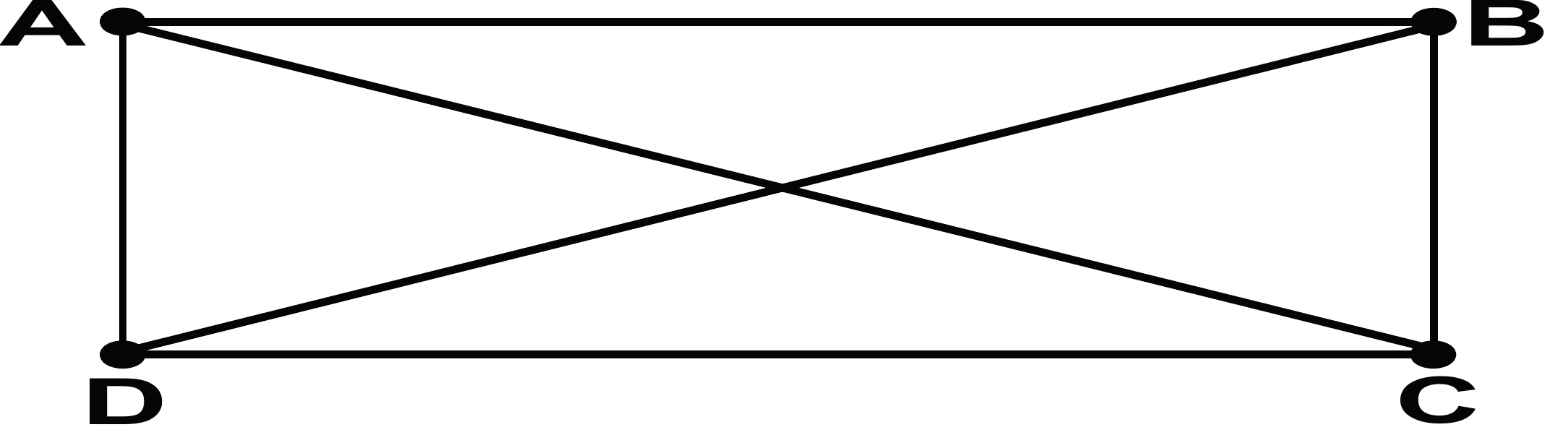}
	\caption{The graph $WH_4 = \textbf{ABCD}$.}
	\label{fig5B}
\end{figure}

\begin{figure}[h]
    \centering
    \includegraphics[width = 0.4\textwidth]{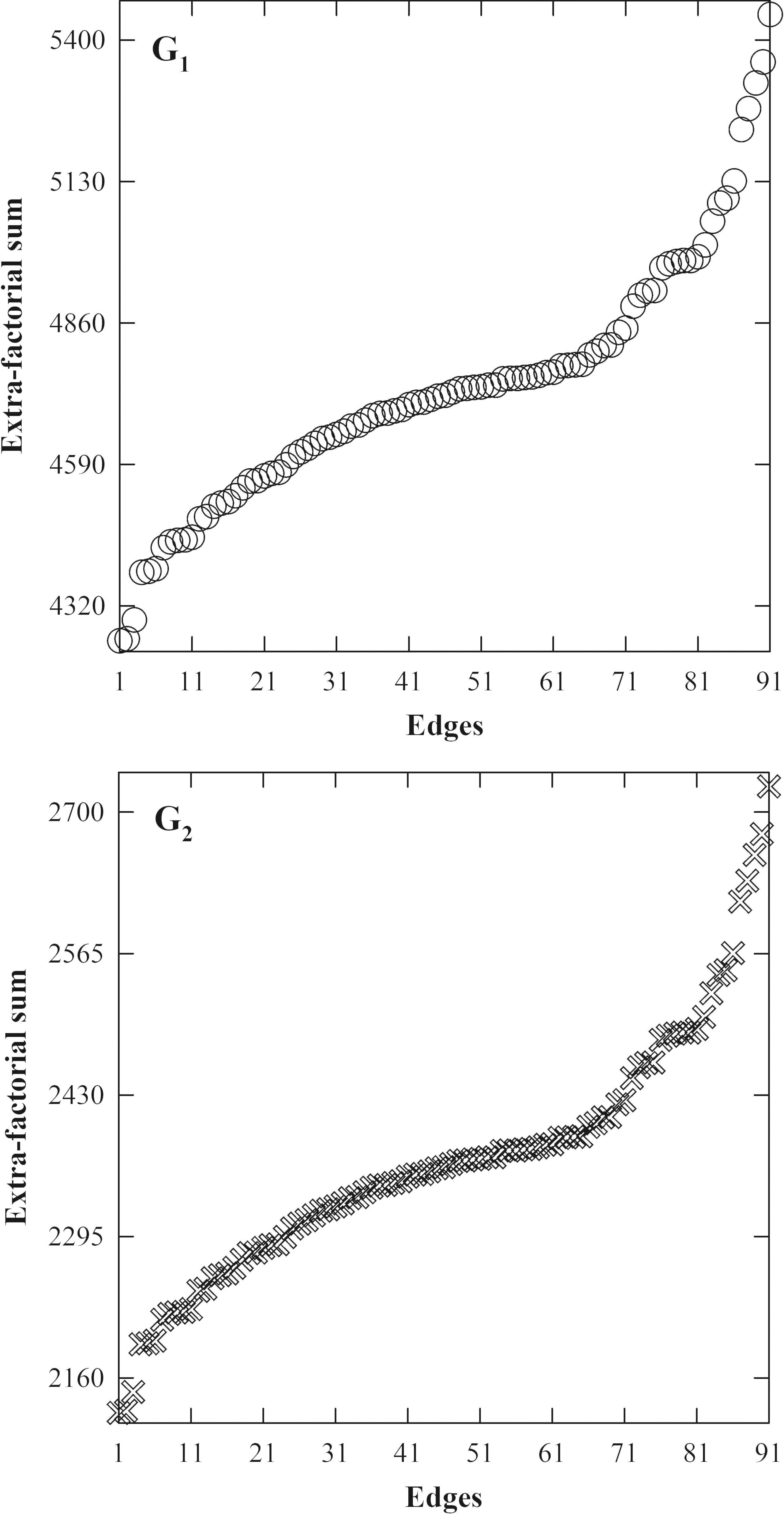}
	\caption{Ranked extra-factorial sum values for two $WH_{14}$ ($G_1$ and $G_2$).}
	\label{fig12}
\end{figure}

For different $WH_n$, the extra-factorial sum can be visualized in a Cartesian chart. Two $WH_{14}$ ($G_1$ and $G_2$) are used as examples: $G_1$ has random edge weights and $G_2$ is a copy of $G_1$ with weights multiplied by $1/2$. The charts of the ranked extra-factorial sums of the 91 edges of each graph are shown in Fig. \ref{fig12}. The first edge has the smallest extra-factorial sum value, while the last (91st) has the largest. This demonstrates that even if the edges have different extra-factorial sums, the corresponding curves are identical. This is because $G_2$ originates from $G_1$ meaning that the graphs are linearly dependent. The similarity in the curves signifies that the cycle length distributions (corresponding to each edge) are also identical. Since the extra-factorial sum provides an overview of the subset of cycles that traverse each edge, future work can focus on absolute or relative similarity \cite{RefJ3} for two or more complete weighted graphs. Consider the graph $WH_4 = \textbf{ABCD}$ (based on Fig. \ref{fig5} but with no vertex \textbf{Y} and edges \textbf{AY, XY, BY, CY}, and where \textbf{A} corresponds to \textbf{A, X} to \textbf{B, B} to \textbf{D, C} to \textbf{C}) (Fig. \ref{fig5B}).

\noindent
The weights are: $w(e_q^{AB}) = 12$, $w(e_q^{AD}) = 7$, $w(e_q^{AC}) = 8$, $w(e_q^{BC}) = 4$, $w(e_q^{BD}) = 5$, and $w(e_q^{CD}) = 2$. The graph has $T_1^4 = \textbf{ACBDA}$ with $l_1 = 24$, $T_2^4 = \textbf{ABCDA}$ with $l_2 = 25$, and $T_3^4 = \textbf{ABDCA}$ with $l_3 = 27$. The arithmetic means of $l_i$ of the $T_i^4$ that traverse each edge are:

\begin{equation}
\frac{(!)\sum_{i=1}^{(n-2)!} l_{i}^{AB}}{(n-2)} = \frac{\sum_{i=1}^{(n-2)!} l_{i}^{AB}}{(n-2)!} = \frac{25+27}{2} = 26
\end{equation}

\begin{equation}
\begin{split}
& \frac{\sum_{i=1}^{(n-2)!} l_{i}^{AC}}{(n-2)!} = 25.5, \frac{\sum_{i=1}^{(n-2)!} l_{i}^{AD}}{(n-2)!} = 24.5 \\
& \frac{\sum_{i=1}^{(n-2)!} l_{i}^{BC}}{(n-2)!} = 24.5, \frac{\sum_{i=1}^{(n-2)!} l_{i}^{BD}}{(n-2)!} = 25.5 \\
& \frac{\sum_{i=1}^{(n-2)!} l_{i}^{CD}}{(n-2)!} = 26
\end{split}
\end{equation}

\noindent
If the weight of each edge (IJ) in $WH_4$ is multiplied by its own value $\frac{\sum_{i=1}^{(n-2)!} l_{i}^{IJ}}{(n-2)!}$, then a new graph $WH_4^{'}$ results. According to Lemma 3, every edge appears $(n - 2)!$ times in the $T_i^4$ of $WH_4^{'}$. Therefore the sum of $l_i$ of the $T_i^4$ ($WH_4^{'}$) is given by:

\begin{equation}
\begin{split}
& (n-2)!e_1^{AB}\frac{\sum_{i=1}^{(n-2)!} l_{i}^{AB}}{(n-2)!}+ \\
& + (n-2)!e_2^{AC}\frac{\sum_{i=1}^{(n-2)!} l_{i}^{AC}}{(n-2)!}+ \\
& + (n-2)!e_3^{AD}\frac{\sum_{i=1}^{(n-2)!} l_{i}^{AD}}{(n-2)!}+ \\
& + (n-2)!e_4^{BC}\frac{\sum_{i=1}^{(n-2)!} l_{i}^{BC}}{(n-2)!}+ \\
& + (n-2)!e_5^{BD}\frac{\sum_{i=1}^{(n-2)!} l_{i}^{BD}}{(n-2)!}+ \\
& + (n-2)!e_6^{CD}\frac{\sum_{i=1}^{(n-2)!} l_{i}^{CD}}{(n-2)!} = \\
& = e_1^{AB}\sum_{i=1}^{(n-2)!} l_{i}^{AB}+e_2^{AC}\sum_{i=1}^{(n-2)!} l_{i}^{AC} +e_3^{AD}\sum_{i=1}^{(n-2)!} l_{i}^{AD}+ \\
& + e_4^{BC}\sum_{i=1}^{(n-2)!} l_{i}^{BC} +e_5^{BD}\sum_{i=1}^{(n-2)!} l_{i}^{BD} +e_6^{CD}\sum_{i=1}^{(n-2)!} l_{i}^{CD}
\end{split}
\end{equation}

\noindent
Substituting the corresponding values, the sum becomes: $12 (25 + 27) + 8 (24 + 27) + 7 (24 + 25) + 4 (24 + 25) + 5 (24 + 27) + 2 (25 + 27)$. This can be written as: $24 (7 + 8 + 5 + 4) + 25 (7 + 4 + 2 + 12) + 27 (8 + 5 + 2 + 12)$ which is: $(576) + (625) + (729)$. Hence, the sum of $l_i$ of the $T_i^4$ ($WH_4^{'}$) is equal to the sum of squared lengths of $T_i^4$ of the initial $WH_4$.

\textbf{Lemma 7}. \textit{The arithmetic mean of $l_i$ of the $T_i^n$ ($WH_n$ or $WH_n^{'}$) is given by:}

\begin{equation}
\begin{split}
& \frac{l_1 + l_2 +l_3 + \cdots + l_{(n-1)!/2}}{(n-1)!/2} = \\ 
& = n\frac{w(e_1)+w(e_2)+\cdots+w(e_{n(n-1)/2})}{n(n-1)/2}
\end{split}
\end{equation}

\textbf{Proof}. According to Lemma 3, the following can be deduced:

\begin{equation}
\begin{split}
& \frac{l_1 + l_2 +l_3 + \cdots + l_{(n-1)!/2}}{(n-1)!/2} = \\
& = \frac{(n-2)![w(e_1)+w(e_2)+\cdots+w(e_{n(n-1)/2})]}{(n-1)!/2}= \\
& = \frac{(n-2)!}{(n-1)!/2}[w(e_1)+w(e_2)+\cdots+w(e_{n(n-1)/2})]= \\ 
& = \frac{n}{n(n-1)/2}[w(e_1)+w(e_2)+\cdots+w(e_{n(n-1)/2})]
\end{split}
\end{equation}

\noindent
This concludes the proof. $\blacksquare$

In a similar framework, the arithmetic mean of the sum of $l_i$ of the $T_i^n$ that do not traverse $e_q$ can be obtained as well. It is interesting to note that these sums can be further analyzed in the context of a controlled change in edge weights, with the aim of analyzing and comparing two or more $WH_n$ using the extra-factorial sum curves. The process of computing the arithmetic mean of the sum of squared lengths of $T_i^n$ ($WH_n$) in combination with the inverse process (given $WH_n^{'}$ to produce an exact corresponding $WH_n$ in which the sum of squared lengths of $T_i^n$ is equal to the sum of $l_i$ of the $T_i^n$ of the initial $WH_n^{'}$) can provide solutions pertinent to the following problem: the existence or not of at least one $T_i^n$ (edge weights having positive and negative values) with negative length. The problem of creating a $WH_n$ from an initial $WH_n^{'}$ is rather complex due to the fact that for each $WH_n^{'}$ there may correspond more than one $WH_n$. Furthermore, an interesting application pertinent to signal processing \cite{RefJ4} is to study the behavior of the extra-factorial sum curves ($WH_n$) that have irrational edge weights (specifically trigonometric numbers) which are functions of time. Yet, another application is in Hopfield neural networks \cite{RefJ5}; these are modelled as $WH_n$ and the extra-factorial sum may contribute to new methods of training the network, e.g. changing the weights in each iteration can incorporate the magnitude of the extra-factorial sum.

\end{document}